\documentclass[12pt,fceqn]{article}
\usepackage{amsmath}
\usepackage{amsfonts}
\usepackage{cite}
\usepackage{amssymb,graphics,graphicx,subfigure,caption2,rotating}
\usepackage{amsmath, latexsym}
\usepackage[amsmath,thmmarks]{ntheorem}
\numberwithin{equation}{section}
\newtheorem{theorem}{Theorem}[section]
\newtheorem{lemma}{Lemma}[section]

\newcommand{\beq}{\begin{equation}}
\newcommand{\eeq}{\end{equation}}
\newcommand{\beqn}{\begin{eqnarray}}
\newcommand{\eeqn}{\end{eqnarray}}
\date{}

\begin{document}

\date{}
\title{Sufficient conditions for the unique solution of a class of new Sylvester-like
absolute value equation\thanks{This research was supported by
National Natural Science Foundation of China (No. 11961082).}}
\author{Shi-Liang Wu\thanks{wushiliang1999@126.com}, Cui-Xia Li\thanks{lixiatkynu@126.com}\\
{\small{\it $^{\dagger}$School of Mathematics, Yunnan Normal
University,}}\\
{\small{\it  Kunming, 650500, PR China }}
}

 \maketitle
\begin{abstract}
%In this paper, a class of new Sylvester-like absolute value equation
%(AVE) $AXB-|CXD|=E$ with $A,C\in \mathbb{R}^{m\times n}$, $B,D\in
%\mathbb{R}^{p\times q}$ and $E\in \mathbb{R}^{m\times q}$ is
%considered, which is quite distinct from the published work by
%Hashemi [Applied Mathematics Letters, 112 (2021) 106818]. Some
%sufficient conditions for the unique solution of the Sylvester-like
%AVE are obtained.

In this paper, a class of new Sylvester-like absolute value equation
(AVE) $AXB-|CXD|=E$ with $A,C\in \mathbb{R}^{m\times n}$, $B,D\in
\mathbb{R}^{p\times q}$ and $E\in \mathbb{R}^{m\times q}$ is
considered, which is quite distinct from the published work by
Hashemi [Applied Mathematics Letters, 112 (2021) 106818]. Some
sufficient conditions for the unique solution of the Sylvester-like
AVE are obtained.


\end{abstract}

\section{Introduction}
As is known, the standard absolute value equation (AVE)
\begin{equation}\label{eq:1}
Ax+|x|=f
\end{equation}
and its general version
\begin{equation}\label{eq:2}
Ax+C|x|=f
\end{equation}
are very strong tools in the field of optimization, including the
complementarity problem, linear programming and convex quadratic
programming, where $A$ and $C$ may be rectangular matrices of the
same order. Based on this, the AVE (\ref{eq:1})/(\ref{eq:2}) has
caused wide public concern over the recent years.

The AVE (\ref{eq:2}) was first introduced in \cite{Rohn} by Rohn.
Therewith, its main research contents consist of two aspects: one is
to multifarious numerical methods for obtaining its numerical
solution, see \cite{Caccetta, Mangasarian, Rohn2, Salkuyeh,Wu, Li,
Guo,Lian,Li2,Li3}, and the other is theoretical analysis, including
the existence of solvability,  bounds for the solutions, various
equivalent reformulations, and so on, see \cite{Rohn3,
Wu2,Mangasarian2,Wu3,Prokopyev,Hlad,Wu4}.

Recently, in \cite{Hashemi}, Hashemi generalized the concept of
absolute value equation and considered the following Sylvester-like
absolute value equation (AVE)
\begin{equation}\label{eq:3}
AXB+C|X|D=F,
\end{equation}
where $A,C\in \mathbb{R}^{m\times n}$, $B,D\in \mathbb{R}^{p\times
q}$, $F\in \mathbb{R}^{m\times q}$ are given. For the Sylvester-like
AVE (\ref{eq:3}), Hashemi in \cite{Hashemi} established some
sufficient conditions for its unique solution. Further, in
\cite{Wang}, Wang and Li considered the Sylvester-like AVE
(\ref{eq:3}) with square coefficient matrices. Some new sufficient
conditions were gained in \cite{Wang}, which are different from the
results in \cite{Hashemi}.

Further, in \cite{Wu5}, Wu found a type of new generalized absolute
value equation (NGAVE), which is below
\begin{equation}\label{eq:4}
Ax-|Bx|=d,
\end{equation}
with $A, B\in \mathbb{R}^{n\times n}$ and $d\in \mathbb{R}^{n}$.
Likewise, Wu in \cite{Wu5} established some necessary and sufficient
conditions  for the unique solution of the NGAVE (\ref{eq:4}).
Clearly, the NGAVE (\ref{eq:4}) is quite different from the AVE
(\ref{eq:2}).

Inspired by the work in \cite{Wu5}, together with the Sylvester-like
AVE (\ref{eq:3}), in this paper, we consider a type of new
Sylvester-like absolute value equation (AVE) below
\begin{equation}\label{eq:5}
AXB-|CXD|=F,
\end{equation}
where $A,C\in \mathbb{R}^{m\times n}$, $B,D\in \mathbb{R}^{p\times
q}$, $F\in \mathbb{R}^{m\times q}$ are given. Here, the new
Sylvester-like AVE (\ref{eq:5}) not only is the generalization form
of the GAVE (\ref{eq:4}), but also is from other fields, such as
interval matrix equations \cite{Neumaier,Seif,
Hashemi2,Hashemi3,Dehghani,Dehghani2}, robust control
\cite{Shashikhin}, and so on. Similar to the Sylvester-like AVE
(\ref{eq:3}), the theory and practice of the new Sylvester-like AVE
(\ref{eq:5}) is still interested and challenged because of the
nonlinear and nondifferentiable term $|CXD|$ in (\ref{eq:5}).  This
is our motivation for this paper. At present, to our knowledge, for
the unique solution of the new Sylvester-like AVE (\ref{eq:5}),  the
necessary and sufficient condition is \emph{vacant}. Based on this,
the goal of the present paper is to fill in  this vacant, gain some
sufficient conditions for the unique solution of the new
Sylvester-like AVE (\ref{eq:5}). What's more, some useful necessary
and sufficient conditions for the unique solution of the new
Sylvester-like AVE (\ref{eq:5}) are obtained with square coefficient
matrices.

\section{Main result}
In this section, we will present some conditions for the unique
solution of the new Sylvester-like AVE (\ref{eq:5}). To achieve this
goal, by using the Kronecker product and the vec operator, the new
Sylvester-like AVE (\ref{eq:5}) can be expressed as the  NGAVE below
\begin{equation}\label{eq:21}
Sx-|Tx|=f
\end{equation}
with $S=B^{T}\otimes A$, $T=D^{T}\otimes C$, $x=vec(X)$ and
$f=vec(F)$, where  `$\otimes$', `$vec$' stand for the Kronecker
product and the vec operator, respectively.

To discuss the sufficient condition for the unique solution of the
new Sylvester-like AVE (\ref{eq:5}), Lemmas 2.1, 2.2, 2.3  and 2.4
are required.
\begin{lemma} \emph{\cite{Wu5}}
Let matrix $A$ in $(\ref{eq:4})$ be nonsingular. If
\begin{equation}\label{eq:22}
\rho((I-2D)BA^{-1})<1
\end{equation}
for any diagonal matrix $D = \mbox{diag}(d_{i})$ with $d_{i}\in [0,
1]$, then the NGAVE $(\ref{eq:4})$ for any $d\in \mathbb{R}^{n}$ has
a unique solution.
\end{lemma}

\begin{lemma}\emph{\cite{Wu5}}
The NGAVE $(\ref{eq:4})$ for any $d\in \mathbb{R}^{n}$ has a unique
solution  if and only if matrix $A+(I-2D)B$ is nonsingular for any
diagonal matrix  $D = \mbox{diag}(d_{i})$ with $d_{i}\in [0, 1]$.
\end{lemma}

\begin{lemma} \emph{\cite{Zhang}}
Let $A,B\in \mathbb{R}^{n\times n}$. Then
\[\sigma_{i}(A+B)\geq\sigma_{i}(A)-\sigma_{1}(B), i=1,2\ldots,n,
\]
where $\sigma_{1}\geq\ldots\geq\sigma_{n}(\geq 0)$ are the singular
values of matrix.
\end{lemma}

Based on Lemmas 2.2 and 2.3, Lemma 2.4 can be obtained.

\begin{lemma}
If
\begin{equation}\label{eq:22}
\sigma_{1}(B)<\sigma_{n}(A)
\end{equation}
where $\sigma_{1}$ and $\sigma_{n}$ denote the largest and smallest
singular value, respectively, then the NGAVE $(\ref{eq:4})$ for any
$d\in \mathbb{R}^{n}$ has a unique solution.
\end{lemma}
\textbf{Proof.} Based on Lemma 2.3, the NGAVE $(\ref{eq:4})$ has a
unique solution for any $d\in \mathbb{R}^{n}$ when the matrix
$A+(I-2D)B$ is nonsingular for any diagonal matrix
$D=\mbox{diag}(d_{i})$ with $0 \leq d_{i}\leq1$. So, let
$\sigma_{n}(A+(I-2D)B)$ stand for the minimal singular value of the
matrix $A+(I-2D)B$. Based on Lemma 2.3, we have
\[
\sigma_{n}(A+(I-2D)B)\geq\sigma_{n}(A)-2\sigma_{1}((I-2D)B).
\]
Since
$\sigma_{1}((I-2D)B)\leq\sigma_{1}((I-2D))\sigma_{1}(B)\leq\sigma_{1}(B)$,
the result in Lemma 2.4 holds under the condition (\ref{eq:22}).
$\hfill{} \Box$

Based on the above lemmas, we can present some conditions for the
unique solution of the new Sylvester-like AVE $(\ref{eq:5})$ for any
$F$. First, by using Lemma 2.1 indirectly, we can obtain the
following result, see Theorem 2.1.

\begin{theorem}
Let $A,B$ be square nonsingular matrix in $(\ref{eq:5})$. If
\begin{equation}\label{eq:23}
\rho((I-2\Lambda)((B^{-1}D)^{T}\otimes CA^{-1}))<1
\end{equation}
for any diagonal matrix  $\Lambda= \mbox{diag}(\lambda_{i})$ with
$\lambda_{i}\in [0, 1]$, then the new Sylvester-like AVE
$(\ref{eq:5})$ has a unique solution for any $F$.
\end{theorem}
\textbf{Proof.} Since
\[
S^{-1}=(B^{T}\otimes A)^{-1}=B^{-T}\otimes A^{-1},
\]
we have
\[
TS^{-1}=(D^{T}\otimes C)(B^{-T}\otimes A^{-1})=D^{T}B^{-T}\otimes
CA^{-1}=(B^{-1}D)^{T}\otimes CA^{-1}.
\]
Based on Lemma 2.1, clearly,  the result in Theorem 2.1 is right.
$\hfill{} \Box$

Clearly, we can use $\rho(((B^{-1}D)^{T}\otimes
CA^{-1})(I-2\Lambda))<1$ instead of the condition (\ref{eq:23})
Theorem 2.1.

In fact, the condition (\ref{eq:23}) in Theorem 2.1 is not easy to
implement in practice. Even if the condition (\ref{eq:23}) in
Theorem 2.1 can be executed, the number of arithmetic operations is
required to compute the spectral radius of the huge matrix. For
instance, for  $A, C \in \mathbb{R}^{m \times m}$ and $B, D \in
\mathbb{R}^{n \times n}$,  the number of arithmetic operations is
$\mathcal{O}(n^{3}m^{3})$. Therefore, the sum of the powers here is
$3+3 = 6$, i.e., the complexity here is sextic. Faced with this
situation, we have to get some conditions that can be detected. By
the simple calculation, we have
\begin{align*}
(I-2\Lambda)((B^{-1}D)^{T}\otimes CA^{-1})&\leq
|(I-2\Lambda)((B^{-1}D)^{T}\otimes CA^{-1}|\\
&\leq|I-2\Lambda||(B^{-1}D)^{T}\otimes CA^{-1})|\\
&\leq|(B^{-1}D)^{T}\otimes CA^{-1})|\\
&=|(B^{-1}D)^{T}|\otimes| CA^{-1}|.
\end{align*}
Note that
\[
\rho(|(B^{-1}D)^{T}|\otimes| CA^{-1}|)=\rho(|(B^{-1}D)^{T}|)\rho(|
CA^{-1}|)=\rho(|B^{-1}D|)\rho(| CA^{-1}|).
\]
So, we have the following result, see Theorem 2.2.
\begin{theorem}
Let $A,B$ be square nonsingular matrix in $(\ref{eq:5})$. If
\begin{equation}\label{eq:24}
\rho(|B^{-1}D|)\rho(| CA^{-1}|)<1,
\end{equation}
then the new Sylvester-like AVE $(\ref{eq:5})$ has a unique solution
for any $F$.
\end{theorem}

In addition, based on Theorem 5.6.10 in \cite{Horn}, we  also have
\begin{align*}
\rho((I-2\Lambda)((B^{-1}D)^{T}\otimes
CA^{-1}))&\leq\sigma_{1}((I-2\Lambda)((B^{-1}D)^{T}\otimes CA^{-1}))\\
&\leq\sigma_{1}(I-2\Lambda)\sigma_{1}((B^{-1}D)^{T}\otimes CA^{-1})\\
&\leq\sigma_{1}((B^{-1}D)^{T}\otimes CA^{-1})\\
&=\sigma_{1}((B^{-1}D)^{T})\sigma_{1}(CA^{-1})\\
&=\sigma_{1}(B^{-1}D)\sigma_{1}(CA^{-1}).
\end{align*}
Based on this, Theorem 2.3 can be obtained.

\begin{theorem}
Let $A,B$ be square nonsingular matrix in $(\ref{eq:5})$. If
\begin{equation}\label{eq:24}
\sigma_{1}(B^{-1}D)\sigma_{1}(CA^{-1})<1,
\end{equation}
then the new Sylvester-like AVE $(\ref{eq:1})$ has a unique solution
for any $F$.
\end{theorem}

Compared with Theorem 2.1, indeed, the conditions of Theorems 2.2
and 2.3 can be easy to execute. Here, it is noted that the
conditions of Theorems 2.2 and 2.3 only work if $A$ and $B$ are
nonsingular. In addition, for a general square matrix $H$, there is
no relations between $\sigma_{1}(H)$ and $\rho(|H|)$ unless one adds
yet more additional requirements. Based on this fact,  Theorem 2.3
sometimes performs better than Theorem 2.2, vice versa.

To obtain the sufficient condition that is more general, based on
Lemma 2.4, Theorem 2.4 can be obtained and its proof is omitted.

\begin{theorem}
If
\begin{equation}\label{eq:24}
\sigma_{1}(C)\sigma_{1}(D)<\sigma_{n}(A)\sigma_{n}(B)
\end{equation}
then the new Sylvester-like AVE $(\ref{eq:1})$ has a unique solution
for any $F$.
\end{theorem}

Compared with Theorems 2.2 and 2.3, Theorem 2.4 not only is fit for
the square matrix, but also is fit for the rectangular matrix. This
implies that the condition (\ref{eq:24}) in Theorem 2.4 is indeed
more general.

It is not difficult to find that all the conditions of Theorems 2.2,
2.3 and 2.4 can be checked. Not only that,  the computational
complexity of all the conditions in Theorems 2.2, 2.3 and 2.4 is
cubic.

By the way, for  $m=n=p=q$ in (\ref{eq:5}), combining Theorems 3.1
and 3.2 in \cite{Wu5}  with Lemma 2.4, we can obtain the following
necessary and sufficient conditions for the unique solution of the
new Sylvester-like AVE $(\ref{eq:5})$, see Theorem 2.5.

\begin{theorem} Let $S=B^{T}\otimes A$, $T=D^{T}\otimes C$. Then
the following statements are equivalent:
\begin{description}
\item  $1.$ the new Sylvester-like AVE $(\ref{eq:1})$ has a unique solution for
any $F\in \mathbb{R}^{n\times n}$;
\item  $2.$  $\{S+T, S-T\}$ has the row $\mathcal{W}$-property;
\item  $3.$  $det(F_{1}(S+T)+F_{2}(S- T))\neq0$ for arbitrary nonnegative diagonal
matrices $F_{1},F_{2} \in \mathbb{R}^{n\times n}$ with
$\mbox{diag}(F_{1}+F_{2}) > 0$;
\item $4$. matrix $(S-T)(S+T)^{-1}$ is a
$P$-matrix (all its principal minors are positive), where matrix
$S+T$ is invertible;
\item $5$. matrix $S+(I-2\Lambda) T$ is nonsingular for any diagonal matrix $\Lambda= \mbox{diag}(\lambda_{i})$ with $\lambda_{i}
\in[0,1]$.
\end{description}
\end{theorem}

Since the order of the matrices $S$ and $T$ in Theorem 2.5 is
$n^{2}\times n^{2}$, the number of arithmetic operations required to
check both parts 2, 3, 4 and 5 of Theorem 2.5 is at-least sextic,
i.e., $\mathcal{O}(n^{6})$. Therefore, the computational complexity
of all the conditions in Theorem 2.5 is at-least
$\mathcal{O}(n^{6})$.

\section{Conclusions}

In this paper, the unique solution of a type of new Sylvester-like
absolute value equation (AVE) $AXB-|CXD|=E$ with $A,C\in
\mathbb{R}^{m\times n}$, $B,D\in \mathbb{R}^{p\times q}$ and $E\in
\mathbb{R}^{m\times q}$ has been discussed. Some useful sufficient
conditions for the unique solution of the new Sylvester-like AVE are
obtained. Particularly, all the sufficient conditions of Theorems
2.2, 2.3 and  2.4 can be checked with a cubic complexity  in the
light of the order of the input matrices.

%\section*{Acknowledgements}
%The authors would like to thank two anonymous referees for providing
%helpful suggestions, which greatly improved the paper.


{\footnotesize
\begin{thebibliography}{99}


\bibitem{Rohn}J. Rohn, Systems of linear interval equations, Linear Algebra
Appl.,
126 (1989) 39-78.


\bibitem{Caccetta} L. Caccetta, B. Qu, G.-L. Zhou, A globally and quadratically
convergent method for absolute value equations, Comput. Optim.
Appl., 48 (2011) 45-58.


\bibitem{Mangasarian}O.L. Mangasarian,  A generalized Newton method for absolute value
equations, Optim. Lett., 3 (2009) 101-108.

\bibitem{Rohn2}J. Rohn, An algorithm for solving the absolute value
equations, Electron. J. Linear Algebra., 18 (2009) 589-599.


\bibitem{Salkuyeh} D.K. Salkuyeh,  The Picard-HSS iteration method for absolute value
equations, Optim. Lett., 8 (2014) 2191-2202.



\bibitem{Wu} S.-L. Wu, C.-X. Li, A special Shift-splitting iterative
method for the absolute value equations, AIMS Math., 5 (2020)
5171-5183.


\bibitem{Li} C.-X. Li, S.-L. Wu, Modified SOR-like iteration method for
absolute value equations, Math. Probl. Eng., 2020 (2020) 9231639.


\bibitem{Guo} P. Guo, S.-L. Wu, C.-X. Li, On SOR-like iteration method for
solving absolute value equations, Appl. Math. Lett., 97 (2019)
107-113.

\bibitem{Lian} Y.-Y. Lian, C.-X. Li, S.-L. Wu, Weaker convergent
results of the generalized Newton method for the generalized
absolute value equations, J. Comput. Appl. Math., 338 (2018)
221-226.


\bibitem{Li2} C.-X. Li, A preconditioned AOR iterative method for the
absolute value equations, Inter. J. Comput. Meth., 14 (2017)
1750016.

\bibitem{Li3} C.-X. Li, A modified generalized Newton method for the
absolute value equations, J. Optim. Theory Appl., 170 (2016)
1055-1059.


\bibitem{Rohn3}J. Rohn, A theorem of the alternatives for the equation
$Ax+B|x|=b$, Linear Multilinear A., 52 (2004) 421-426.

\bibitem{Wu2}S.-L. Wu, C.-X. Li, A note on unique solvability of the absolute value equation, Optim.
Lett., 14 (2020) 1957-1960.

\bibitem{Mangasarian2} O.L. Mangasarian,  R.R. Meyer, Absolute value equations, Linear
Algebra Appl.,  419 (2006) 359-367.


\bibitem{Wu3}S.-L. Wu, C.-X. Li, The unique solution of the absolute value
equations, Appl. Math. Lett., 76 (2018) 195-200.

\bibitem{Prokopyev}O. Prokopyev, On equivalent reformulations for absolute value
equations, Comput. Optim. Appl., 44 (2009) 363-372.


\bibitem{Hlad}M. Hlad\'{\i}k, Bounds for the solutions of absolute value
equations, Comput. Optim. Appl., 69 (2018) 243-266.

\bibitem{Wu4} S.-L. Wu, S.-Q. Shen, On the unique solution of the generalized
absolute value equation, Optim. Lett., 2020,
https://doi.org/10.1007/s11590-020-01672-2.

\bibitem{Hashemi}B. Hashemi, Sufficient conditions for the solvability of a
Sylvester-like absolute value matrix equation, Appl. Math. Lett.,
112 (2021) 106818.

\bibitem{Wang}L.-M. Wang, C.-X. Li, New sufficient conditions for the unique solution of a square
Sylvester-like absolute value equation, Appl. Math. Lett.,  116
(2021) 106966.


\bibitem{Wu5}S.-L. Wu, The unique solution of a class of the new generalized absolute value
equation, Appl. Math. Lett., 116 (2021) 107029.

\bibitem{Neumaier} A. Neumaier, Interval Methods for Systems of Equations, Cambridge
University Press, Cambridge, 1990.

\bibitem{Seif} N.P. Seif, S.A. Hussein, A.S. Deif, The interval Sylvester
equation, Computing., 52 (1994) 233-244.

\bibitem{Hashemi2} B. Hashemi, M. Dehghan, Results concerning interval linear
systems with multiple right-hand sides and the interval matrix
equation $AX = B$, J. Comput. Appl. Math., 235 (2011) 2969-2978.

\bibitem{Hashemi3} B. Hashemi, M. Dehghan, The interval Lyapunov matrix equation:
analytical results and an efficient numerical technique for outer
estimation of the united solution set, Math. Comput. Model., 55
(2012) 622-633.

\bibitem{Dehghani} M. Dehghani-Madiseh, M. Hlad\'{\i}k,  Efficient approaches for enclosing the united
solution set of the interval generalized Sylvester matrix equations,
Appl. Numer. Math., 126 (2018) 18-33.

\bibitem{Dehghani2}M. Dehghani-Madiseh, M. Dehghan, Generalized solution sets of
the interval generalized Sylvester matrix equation $\sum^{p}_{i =1}
A_{i}X_{i} +\sum^{q}_{j =1}Y_{j}B_{j }= C$ and some approaches for
inner and outer estimations, Comput.  Math. Appl., 68 (2014)
1758-1774.


\bibitem{Shashikhin} V.N. Shashikhin, Robust assignment of poles in large-scale
interval systems, Autom. Rem. Contr., 63 (2002) 200-208.


\bibitem{Zhang}F.-Z. Zhang, Matrix Theory: Basic results and techniques (Second edition).
Springer, New York, 2011.


\bibitem{Horn} R.A. Horn,  C.R. Johnson, Matrix Analysis. Cambridge University Press,
Cambridge, 1986.

\end{thebibliography}}
\end{document}